# Nonparametric Bayesian model selection and averaging

## Subhashis Ghosal

*Department of Statistics, North Carolina State University, 2501 Founders Drive, Raleigh,
North Carolina 27695-8203, USA, e-mail:* `ghosal@stat.ncsu.edu`

## Jüri Lember[*]

*Institute of Mathematical Statistics, University of Tartu, J. Liivi 2, 50409 Tartu, Estonia,
e-mail:* `juri.lember@ut.ee`

## Aad van der Vaart

*Department of Mathematics, Faculty of Sciences, Vrije Universiteit, De Boelelaan 1081a,
1081 HV Amsterdam, Netherlands, e-mail:* `aad@cs.vu.nl`

**Abstract:** We consider nonparametric Bayesian estimation of a probability density $p$ based on a random sample of size $n$ from this density using a hierarchical prior. The prior consists, for instance, of prior weights on the regularity of the unknown density combined with priors that are appropriate given that the density has this regularity. More generally, the hierarchy consists of prior weights on an abstract model index and a prior on a density model for each model index. We present a general theorem on the rate of contraction of the resulting posterior distribution as $n \to \infty$, which gives conditions under which the rate of contraction is the one attached to the model that best approximates the true density of the observations. This shows that, for instance, the posterior distribution can adapt to the smoothness of the underlying density. We also study the posterior distribution of the model index, and find that under the same conditions the posterior distribution gives negligible weight to models that are bigger than the optimal one, and thus selects the optimal model or smaller models that also approximate the true density well. We apply these result to log spline density models, where we show that the prior weights on the regularity index interact with the priors on the models, making the exact rates depend in a complicated way on the priors, but also that the rate is fairly robust to specification of the prior weights.

---

[*]Research carried out at Eurandom, Eindhoven, 2003–2004, and supported by Estonian Science Foundation Grant 7553.

## 1. Introduction

It is well known that the selection of a suitable "bandwidth" is crucial in non-parametric estimation of densities. Within a Bayesian framework it is natural to put a prior on bandwidth and let the data decide on a correct bandwidth through the corresponding posterior distribution. More generally a Bayesian procedure might consist of the specification of a suitable prior on a statistical model that is "correct" if the true density possesses a certainly regularity level, together with the specification of a prior on the regularity. Such a hierarchical Bayesian procedure fits naturally within the framework of *adaptive estimation*, which focuses on constructing estimators that automatically choose a best model from a given set of models. Given a collection of models an estimator is said to be *rate-adaptive* if it attains the rate of convergence that would have been attained had only the best model been used. For instance, the minimax rate of convergence for estimating a density on $[0, 1]^d$ that is known to have $\alpha$ derivatives is $n^{-\alpha/(2\alpha+d)}$. An estimator would be *rate-adaptive* to the set of models consisting of all smooth densities if it attained the rate $n^{-\alpha/(2\alpha+d)}$ whenever the true density is $\alpha$-smooth, for any $\alpha > 0$. (See e.g. Tsybakov [2004].)

In this paper we present a general result on adaptation for density estimation within the Bayesian framework. The observations are a random sample $X_1, \ldots, X_n$ from a density on a given measurable space. Given a countable collection of density models $\mathcal{P}_{n,\alpha}$, indexed by a parameter $\alpha \in A_n$, each provided with a prior distribution $\Pi_{n,\alpha}$, and a prior distribution $\lambda_n$ on $A_n$, we consider the posterior distribution relative to the prior that first chooses $\alpha$ according to $\lambda_n$ and next $p$ according to $\Pi_{n,\alpha}$ for the chosen $\alpha$. The index $\alpha$ may be a regularity parameter, but in the general result it may be arbitrary. Thus the overall prior is a probability measure on the set of probability densities, given by

$$\Pi_n = \sum_{\alpha \in A_n} \lambda_{n,\alpha} \Pi_{n,\alpha}. \tag{1.1}$$

Given this prior distribution, the corresponding posterior distribution is the random measure

$$\begin{aligned} \Pi_n(B|X_1, \ldots, X_n) &= \frac{\int_B \prod_{i=1}^n p(X_i) \, d\Pi_n(p)}{\int \prod_{i=1}^n p(X_i) \, d\Pi_n(p)} \\ &= \frac{\sum_{\alpha \in A_n} \lambda_{n,\alpha} \int_{p \in \mathcal{P}_{n,\alpha}: p \in B} \prod_{i=1}^n p(X_i) \, d\Pi_{n,\alpha}(p)}{\sum_{\alpha \in A_n} \lambda_{n,\alpha} \int_{p \in \mathcal{P}_{n,\alpha}} \prod_{i=1}^n p(X_i) \, d\Pi_{n,\alpha}(p)}. \end{aligned} \tag{1.2}$$

Of course, we make appropriate (measurability) conditions to ensure that this expression is well defined.

We say that the posterior distributions have *rate of convergence at least* $\varepsilon_n$ if, for every sufficiently large constant $M$, as $n \to \infty$, in probability,

$$\Pi_n\big(d(p, p_0) > M \varepsilon_n | X_1, \ldots, X_n\big) \to 0.$$

Here the distribution of the random measure (1.2) is evaluated under the assumption that $X_1, \ldots, X_n$ are an i.i.d. sample from $p_0$, and $d$ is a distance on the set of densities. Throughout the paper this distance is assumed to be bounded above by the Hellinger distance and generate convex balls. (For instance, the Hellinger or $L_1$-distance, or the $L_2$-distance if the densities are uniformly bounded.) Thus we study the asymptotics of the posterior distribution in the frequentist sense.

The aim is to prove a result of the following type. For a given $p_0$ there exists a best model $\mathcal{P}_{n,\beta_n}$ that gives a posterior rate $\varepsilon_{n,\beta_n}$ if it would be combined with the prior $\Pi_{n,\beta_n}$. The hierachical Bayesian procedure would adapt to the set of models if the posterior distributions (1.2), which are based on the mixture prior (1.1), have the same rate of convergence for this $p_0$, for any $p_0$ in some model $\mathcal{P}_{n,\alpha}$. Technically, the first main result is Theorem 2.1 in Section 2.

This sense of Bayesian adaptation refers to the "full" posterior, both its centering and its spread. As noted by Belitser and Ghosal [2003] suitably defined centers of these posteriors would yield adaptive point estimators.

The posterior distribution can be viewed as a mixture of the posterior distributions on the various models, with the weights given by the posterior distribution of the model index. Our second main result, Theorem 3.1 concerns the posterior distribution of the model index. It shows that models that are "bigger" than the optimal model asymptotically achieve zero posterior mass. On the other hand, under our conditions the posterior may distribute its mass over a selection of smaller models, provided that these can approximate the true distribution well.

In the situation that there are precisely two models this phenomenon can be conveniently described by the *Bayes factor* of the two models. We provide simple sufficient conditions for the Bayes factor to select the "true" model with probability tending to one. This consistency property is especially relevant for Bayesian goodness of fit testing against a nonparametric alternative. A computationally advantageous method of such goodness of fit test was developed by Berger and Guglielmi [2001] using a mixture of Polya tree prior on the nonparametric alternative. Asymptotic properties of Bayes factors for nested regular parametric models have been well studied beginning with the pioneering work by Schwarz [1978], who also introduced the Bayesian information criterion. However, large sample properties of Bayes factors when at least one model is infinite dimensional appear to be unknown except in special cases. The paper Dass and Lee [2004] showed consistency of Bayes factors when one of the models is a singleton and the prior for the other model assigns positive probabilities to the Kullback-Leibler neighborhoods of the true density, popularly known as the Kullback-Leibler property. The paper Walker et al. [2004] showed (in particular) that if the prior on one model has the Kullback-Leibler property and the other does not, then the Bayes factor will asymptotically favour the model with the Kullback-Leibler property. Unfortunately, the proof of Dass and Lee [2004] does not generalize to general null models and frequently both priors will have the Kullback-Leibler property, precluding the application of Walker et al. [2004]. In Sections 3 and 4 we study these issues in general.

The present paper is an extension of the paper Ghosal et al. [2003], which studies adaptation to finitely many models of splines with a uniform weight on the models. In the present paper we derive a result for general models, possibly infinitely many, and investigate different model weights. Somewhat surprisingly we find that both weights that give more prior mass to small models and weights that downweight small models may lead to adaptation.

Related work on Bayesian adaptation was carried out by Huang [2004], who considers adaptation using scales of finite-dimensional models, and Lember and van der Vaart (2007), who consider special weights that downweight large models. Our methods of proof borrow from Ghosal et al. [2000].

The paper is organized as follows. After stating the main theorems on adaptation and model selection and some corollaries in Sections 2 and 3, we investigate adaptation in detail in the context of log spline models in Section 5, and we consider the Bayes factors for testing a finite- versus an infinite-dimensional model in detail in Section 4. The proof of the main theorems is given in Section 6, and further technical proofs and complements are given in Section 7.

### 1.1. Notation

Throughout the paper the data are a random sample $X_1, \ldots, X_n$ from a probability measure $P_0$ on a measurable space $(\mathcal{X}, \mathcal{A})$ with density $p_0$ relative to a given reference measure $\mu$ on $(\mathcal{X}, \mathcal{A})$. In general we write $p$ and $P$ for a density and the corresponding probability measure. The Hellinger distance between two densities $p$ and $q$ relative to $\mu$ is defined as $h(p, q) = \|\sqrt{p} - \sqrt{q}\|_2$, for $\|\cdot\|_2$ the norm of $L_2(\mu)$. The *$\varepsilon$-covering numbers* and *$\varepsilon$-packing numbers* of a metric space $(\mathcal{P}, d)$, denoted by $N(\varepsilon, \mathcal{P}, d)$ and $D(\varepsilon, \mathcal{P}, d)$, are defined as the minimal numbers of balls of radius $\varepsilon$ needed to cover $\mathcal{P}$, and the maximal number of $\varepsilon$-separated points, respectively.

For each $n \in \mathbb{N}$ the index set is a countable set $A_n$, and for every $\alpha \in A_n$ the set $\mathcal{P}_{n,\alpha}$ is a set of $\mu$-probability densities on $(\mathcal{X}, \mathcal{A})$ equipped with a $\sigma$-field such that the maps $(x, p) \mapsto p(x)$ are measurable. Furthermore, $\Pi_{n,\alpha}$ denotes a probability measure on $\mathcal{P}_{n,\alpha}$, and $\lambda_n = (\lambda_{n,\alpha} : \alpha \in A_n)$ is a probability measure on $A_n$. We define

$$B_{n,\alpha}(\varepsilon) = \left\{ p \in \mathcal{P}_{n,\alpha} : -P_0 \log \frac{p}{p_0} \leq \varepsilon^2, P_0 \left( \log \frac{p}{p_0} \right)^2 \leq \varepsilon^2 \right\},$$

$$C_{n,\alpha}(\varepsilon) = \left\{ p \in \mathcal{P}_{n,\alpha} : d(p, p_0) \leq \varepsilon \right\}. \tag{1.3}$$

Throughout the paper $\varepsilon_{n,\alpha}$ are given positive numbers with $\varepsilon_{n,\alpha} \to 0$ as $n \to \infty$. These may be thought of as the rate attached to the model $\mathcal{P}_{n,\alpha}$ if this is (approximately) correct.

The notation $a \lesssim b$ means that $a \leq Cb$ for a constant $C$ that is universal or fixed in the proof. For sequences $a_n$ and $b_n$ we write $a_n \ll b_n$ if $a_n/b_n \to 0$ and $a_n \gg 0$ if $a_n > 0$ for every $n$ and $\liminf a_n > 0$. For a measure $P$ and a measurable function $f$ we write $Pf$ for the integral of $f$ relative to $P$.

## 2. Adaptation

For $\beta_n$ a given element of $A_n$, thought to be the index of a best model for a given fixed true density $p_0$, we split the index set in the indices that give a faster or slower rate: for a fixed constant $H \geq 1$,

$$A_{n,\gtrsim\beta_n} := \Big\{ \alpha \in A_n : \varepsilon_{n,\alpha}^2 \leq H\varepsilon_{n,\beta_n}^2 \Big\},$$

$$A_{n,<\beta_n} := \Big\{ \alpha \in A_n : \varepsilon_{n,\alpha}^2 > H\varepsilon_{n,\beta_n}^2 \Big\}.$$

Even though we do not assume that $A_n$ is ordered, we shall write $\alpha \gtrsim \beta_n$ and $\alpha < \beta_n$ if $\alpha$ belongs to the sets $A_{n,\gtrsim\beta_n}$ or $A_{n,<\beta_n}$, respectively. The set $A_{n,\gtrsim\beta_n}$ contains $\beta_n$ and hence is never empty, but the set $A_{n,<\beta_n}$ can be empty (if $\beta_n$ is the "smallest" possible index). In the latter case conditions involving $\alpha < \beta_n$ are understood to be automatically satisfied.

The assumptions of the following theorem are reminiscent of the assumptions in Ghosal et al. [2000], and entail a bound on the complexity of the models and a condition on the concentration of the priors. The complexity bound is exactly as in Ghosal et al. [2000] and takes the form: for some constants $E_\alpha$,

$$\sup_{\varepsilon \geq \varepsilon_{n,\alpha}} \log N\Big(\frac{\varepsilon}{3}, C_{n,\alpha}(2\varepsilon), d\Big) \leq E_\alpha n\varepsilon_{n,\alpha}^2, \qquad \alpha \in A_n. \tag{2.1}$$

The conditions on the priors involve comparisons of the prior masses of balls of various sizes in various models. These conditions are split in conditions on the models that are smaller or bigger than the best model: for given constants $\mu_{n,\alpha}, L, H, I,$

$$\frac{\lambda_{n,\alpha}}{\lambda_{n,\beta_n}} \frac{\Pi_{n,\alpha}\big(C_{n,\alpha}(i\varepsilon_{n,\alpha})\big)}{\Pi_{n,\beta_n}\big(B_{n,\beta_n}(\varepsilon_{n,\beta_n})\big)} \leq \mu_{n,\alpha} e^{Li^2 n\varepsilon_{n,\alpha}^2}, \qquad \alpha < \beta_n, \quad i \geq I, \tag{2.2}$$

$$\frac{\lambda_{n,\alpha}}{\lambda_{n,\beta_n}} \frac{\Pi_{n,\alpha}\big(C_{n,\alpha}(i\varepsilon_{n,\beta_n})\big)}{\Pi_{n,\beta_n}\big(B_{n,\beta_n}(\varepsilon_{n,\beta_n})\big)} \leq \mu_{n,\alpha} e^{Li^2 n\varepsilon_{n,\beta_n}^2}, \qquad \alpha \gtrsim \beta_n, \quad i \geq I. \tag{2.3}$$

A final condition requires that the prior mass in a ball of radius $\varepsilon_{n,\alpha}$ in a big model (i.e. small $\alpha$) is significantly smaller than in a small model: for some constants $I, B,$

$$\sum_{\alpha \in A_n : \alpha < \beta_n} \frac{\lambda_{n,\alpha}}{\lambda_{n,\beta_n}} \frac{\Pi_{n,\alpha}\big(C_{n,\alpha}(IB\varepsilon_{n,\alpha})\big)}{\Pi_{n,\beta_n}\big(B_{n,\beta_n}(\varepsilon_{n,\beta_n})\big)} = o\big(e^{-2n\varepsilon_{n,\beta_n}^2}\big). \tag{2.4}$$

Let $K$ be the universal testing constant. According to the assertion of Lemma 6.2 it can certainly be taken equal to $K = 1/9$.

**Theorem 2.1.** *Assume there exist positive constants $B, E_\alpha, L, H \geq 1, I > 2$ such that (2.1), (2.2), (2.3) and (2.4) hold, and, constants $E$ and $\underline{E}$ such that $E \geq \sup_{\alpha \in A_n : \alpha \gtrsim \beta_n} E_\alpha \varepsilon_{n,\alpha}^2 / \varepsilon_{n,\beta_n}^2$ and $\underline{E} \geq \sup_{\alpha \in A_n : \alpha < \beta_n} E_\alpha$ (with $\underline{E} = 0$ if $A_{n,<\beta_n} = \emptyset$),*

$$B > \sqrt{H}, \quad KB^2 > (H\underline{E}) \vee E + 1, \quad B^2 I^2 (K - 2L) > 3.$$

*Furthermore, assume that $\sum_{\alpha \in A_n} \sqrt{\mu_{n,\alpha}} \le \exp[n\varepsilon_{n,\beta_n}^2]$. If $\beta_n \in A_n$ for every $n$ and satisfies $n\varepsilon_{n,\beta_n}^2 \to \infty$, then the posterior distribution (1.2) satisfies*

$$P_0^n \Pi_n \Big( p \colon d(p, p_0) \ge IB\,\varepsilon_{n,\beta_n} | X_1, \cdots, X_n \Big) \to 0.$$

The proof of the theorem is deferred to Section 6.

In many situations (although not in the main examples of the present paper) relatively crude bounds on the prior mass bounds (2.2), (2.3) and (2.4) are sufficient. In particular, the following lower bound is often useful: for a positive constant $F$,

$$\Pi_{n,\beta_n} \big( B_{n,\beta_n}(\varepsilon_{n,\beta_n}) \big) \ge \exp[-Fn\varepsilon_{n,\beta_n}^2]. \tag{2.5}$$

This correspond to the "crude" prior mass condition of Ghosal et al. [2000]. Combined with the trivial bound 1 on the probabilities $\Pi_{n,\alpha}(C)$ in (2.2) and (2.3), we see that these conditions hold (for sufficiently large $I$) if, for all $\alpha \in A_n$,

$$\frac{\lambda_{n,\alpha}}{\lambda_{n,\beta_n}} \le \mu_{n,\alpha} e^{n(\varepsilon_{n,\alpha}^2 \vee \varepsilon_{n,\beta_n}^2)}. \tag{2.6}$$

This appears to be a mild requirement. On the other hand, the similarly adapted version of condition (2.4) still requires that

$$\sum_{\alpha \in A_n : \alpha < \beta_n} \frac{\lambda_{n,\alpha}}{\lambda_{n,\beta_n}} \Pi_{n,\alpha} \big( C_{n,\alpha}(IB\varepsilon_{n,\alpha}) \big) = o\big( e^{-(F+2)n\varepsilon_{n,\beta_n}^2} \big). \tag{2.7}$$

Such a condition may be satisfied because the prior probabilities $\Pi_{n,\alpha} \big( C_{n,\alpha}(IB\varepsilon_{n,\alpha}) \big)$ are very small. For instance, a reverse bound of the type (2.5) for $\alpha$ instead of $\beta_n$ would yield this type of bound for fairly general model weights $\lambda_{n,\alpha}$, since $\varepsilon_{n,\alpha} \ge H\varepsilon_{n,\beta_n}$ for $\alpha < \beta_n$. Alternatively, the condition could be forced by choice of the model weights $\lambda_{n,\alpha}$, for general priors $\Pi_{n,\alpha}$. For instance, in Section 5.3 we consider weights of the type

$$\lambda_{n,\alpha} = \frac{\mu_\alpha \exp[-Cn\varepsilon_{n,\alpha}^2]}{\sum_\alpha \mu_\alpha \exp[-Cn\varepsilon_{n,\alpha}^2]}. \tag{2.8}$$

Such weights were also considered in Lember and van der Vaart [2007], who discuss several other concrete examples. For reference we codify the preceding discussion as a lemma.

**Lemma 2.1.** *Conditions (2.5), (2.6) and (2.7) are sufficient for (2.2), (2.3) and (2.4).*

Theorem 2.1 excludes the case that $\varepsilon_{n,\beta_n}$ is equal to the "parametric rate" $1/\sqrt{n}$. To cover this case the statement of the theorem must be slightly adapted. The proof of the following theorem is given in Section 6.

**Theorem 2.2.** *Assume there exist positive constants $B, E_\alpha, L < K/2, H \ge 1, I$ such that (2.1), (2.2), (2.3) and (2.4) hold for every sufficiently large $I$.*

*Furthermore, assume that* $\sum_{\alpha \in A_n} \sqrt{\mu_{n,\alpha}} = O(1)$. *If* $\beta_n \in A_n$ *for every* $n$ *and* $\varepsilon_{n,\beta_n} = 1/\sqrt{n}$, *then the posterior distribution (1.2) satisfies, for every* $I_n \to \infty$,

$$P_0^n \Pi_n \Big( p : d(p, p_0) \geq I_n \varepsilon_{n,\beta_n} | X_1, \cdots, X_n \Big) \to 0.$$

For further understanding it is instructive to apply the theorems to the situation of two models, say $\mathcal{P}_{n,1}$ and $\mathcal{P}_{n,2}$ with rates $\varepsilon_{n,1} > \varepsilon_{n,2}$. For simplicity we shall also assume (2.5) and use universal constants.

**Corollary 2.1.** *Assume that (2.1) holds for* $\alpha \in A_n = \{1, 2\}$ *and sequences* $\varepsilon_{n,1} > \varepsilon_{n,2}$.

*(1) If* $\Pi_{n,1} \big( B_{n,1}(\varepsilon_{n,1}) \big) \geq e^{-n\varepsilon_{n,1}^2}$ *and* $\lambda_{n,2}/\lambda_{n,1} \leq e^{n\varepsilon_{n,1}^2}$, *then the posterior rate of contraction is at least* $\varepsilon_{n,1}$.

*(2) If* $\Pi_{n,2} \big( B_{n,2}(\varepsilon_{n,2}) \big) \geq e^{-n\varepsilon_{n,2}^2}$ *and* $\lambda_{n,2}/\lambda_{n,1} \geq e^{-n\varepsilon_{n,1}^2}$ *and, moreover,* $\Pi_{n,1} \big( C_{n,1}(I\varepsilon_{n,1}) \big) \leq (\lambda_{n,2}/\lambda_{n,1}) o(e^{-3n\varepsilon_{n,2}^2})$ *for every* $I$, *then the posterior rate of contraction is at least* $\varepsilon_{n,2}$.

*Proof.* We apply the preceding theorems with $\beta_n = 1$, $A_{n,<\beta_n} = \emptyset$ and $A_{n,\gtrsim\beta_n} = \{1, 2\}$ in case (1) and $\beta_n = 2$, $A_{n,<\beta_n} = \{1\}$ and $A_{n,\gtrsim\beta_n} = \{2\}$ in case (2), both times with $H = 1$ and $\mu_{n,1} = \mu_{n,2} = 1$. ∎

Statement (1) of the corollary gives the slower $\varepsilon_{n,1}$ of the two rates under the assumption that the bigger model satisfies the prior mass condition (2.5) and a condition on the weights $\lambda_{n,i}$ that ensures that the smaller model is not overly downweighted. The latter condition is very mild, as it allows the weights of the two models too be very different. Apart from this, statement (1) is not surprising, and could also be obtained from nonadaptive results on posterior rates of contraction, as in Ghosal et al. [2000].

Statement (2) gives the faster rate of contraction $\varepsilon_{n,2}$ under the condition that the smaller model satisfies the prior mass condition (2.5), an equally mild condition on the relative weights of the two models, and an additional condition on the prior weight $\Pi_{n,1} \big( C_{n,1}(I\varepsilon_{n,1}) \big)$ that the bigger model attaches to neighbourhoods of the true distribution. If this would be of the expected order $\exp(-Fn\varepsilon_{n,1}^2)$, then the conditions on the weights $\lambda_{n,1}$ and $\lambda_{n,2}$ in the union of (1) and (2) can be summarized as

$$e^{-n\varepsilon_{n,1}^2} \leq \frac{\lambda_{n,2}}{\lambda_{n,1}} \leq e^{n\varepsilon_{n,1}^2}.$$

This is a remarkably big range of weights. One might conclude that Bayesian methods are very robust to the prior specification of model weights. One might also more cautiously guess that rate-asymptotics do not yield a complete picture of the performance of the various priors (even though rates are considerably more informative than consistency results).

**Remark 2.1.** *The entropy condition (2.1) can be relaxed to the same condition on a submodel* $\mathcal{P}'_{n,\alpha} \subset \mathcal{P}_{n,\alpha}$ *that carries most of the prior mass, in the sense*

*that*

$$\frac{\sum_\alpha \lambda_{n,\alpha} \Pi_{n,\alpha}(\mathcal{P}_{n,\alpha} - \mathcal{P}'_{n,\alpha})}{\lambda_{n,\beta_n} \Pi_{n,\beta_n}\big(B_{n,\beta_n}(\varepsilon_{n,\beta_n})\big)} = o\big(e^{-2n\varepsilon^2_{n,\beta_n}}\big).$$

*This follows, because in that case the posterior will concentrate on $\cup_\alpha \mathcal{P}'_{n,\alpha}$ (see Ghosal and van der Vaart [2007a], Lemma 1). This relaxation has been found useful in several papers on nonadaptive rates. In the present context, it seems that the condition would only be natural if it is valid for the index $\beta_n$ that gives the slowest of all rates $\varepsilon_{n,\alpha}$.*

## 3. Model Selection

The theorems in the preceding section concern the concentration of the posterior distribution on the set of densities relative to the metric $d$. In this section we consider the posterior distribution of the index parameter $\alpha$, within the same set-up. The proof of the following theorem is given in Section 6.

Somewhat abusing notation (cf. (1.2)), we write, for any set $B \subset A_n$ of indices,

$$\Pi_n(B|X_1, \ldots, X_n) = \frac{\sum_{\alpha \in B} \lambda_{n,\alpha} \int \prod_{i=1}^n p(X_i) \, d\Pi_{n,\alpha}(p)}{\sum_{\alpha \in A_n} \lambda_{n,\alpha} \int \prod_{i=1}^n p(X_i) \, d\Pi_{n,\alpha}(p)}.$$

**Theorem 3.1.** *Under the conditions of Theorem 2.1,*

$$P_0^n \Pi_n\big(A_{n,<\beta_n}|X_1, \cdots, X_n\big) \to 0,$$
$$P_0^n \Pi_n\big(\alpha \in A_{n,\gtrsim\beta_n} : d(p_0, \mathcal{P}_{n,\alpha}) > IB\,\varepsilon_{n,\beta_n}|X_1, \cdots, X_n\big) \to 0.$$

*Under the conditions of Theorem 2.2 this is true with $IB$ replaced by $I_n$, for any $I_n \to \infty$.*

The first assertion of the theorem is pleasing. It can be interpreted in the sense that the models that are bigger than the model $\mathcal{P}_{n,\beta_n}$ that contains the true distribution eventually receive negligible posterior weight. The second assertion makes a similar claim about the smaller models, but it is restricted to the smaller models that keep a certain distance to the true distribution. Such a restriction appears not unnatural, as a small model that can represent the true distribution well ought to be favoured by the posterior: the posterior looks at the data through the likelihood and hence will judge a model by its approximation properties rather than its parametrization. That big models with similarly good approximation properties are not favoured is caused by the fact that (under our conditions) the prior mass on the big models is more spread out, yielding relatively little prior mass near good approximants within the big models.

It is again insightful to specialize the theorem to the case of two models, and simplify the prior mass conditions to (2.5). The behaviour of the posterior of the model index can then be described through the *Bayes factor*

$$\mathrm{BF}_n = \frac{\lambda_{n,2} \int \prod_{i=1}^n p(X_i) \, \Pi_{n,2}(p)}{\lambda_{n,1} \int \prod_{i=1}^n p(X_i) \, \Pi_{n,1}(p)}.$$

**Corollary 3.1.** *Assume that (2.1) holds for $\alpha \in A_n = \{1, 2\}$ and sequences $\varepsilon_{n,1} > \varepsilon_{n,2}$.*

(1) *If $\Pi_{n,1}\big(B_{n,1}(\varepsilon_{n,1})\big) \geq e^{-n\varepsilon_{n,1}^2}$ and $\lambda_{n,2}/\lambda_{n,1} \leq e^{n\varepsilon_{n,1}^2}$ and $d(p_0, \mathcal{P}_{n,2}) \geq I_n \varepsilon_{n,1}$ for every $n$ and some $I_n \to \infty$, then $\mathrm{BF}_n \to 0$ in $P_0^n$-probability.*

(2) *If $\Pi_{n,2}\big(B_{n,2}(\varepsilon_{n,2})\big) \geq e^{-n\varepsilon_{n,2}^2}$ and $\lambda_{n,2}/\lambda_{n,1} \geq e^{-n\varepsilon_{n,1}^2}$ and $\Pi_{n,1}\big(C_{n,1} \times (I\varepsilon_{n,1})\big) \leq (\lambda_{n,2}/\lambda_{n,1}) o(e^{-3n\varepsilon_{n,2}^2})$ for every $I$, then $\mathrm{BF}_n \to \infty$ in $P_0^n$-probability.*

*Proof.* The Bayes factor tends to 0 or $\infty$ if the posterior probability of model $\mathcal{P}_{n,2}$ or $\mathcal{P}_{n,1}$ tends to zero, respectively. Therefore, we can apply Theorem 3.1 with the same choices as in the proof of Corollary 2.1. ∎

In particular, if the two models are equally weighted ($\lambda_{n,1} = \lambda_{n,2}$), the models satisfy (2.1) and the priors satisfy (2.5), then the Bayes factors are asymptotically consistent if

$$d(p_0, \mathcal{P}_{n,1}) \gg \varepsilon_{n,1},$$

$$\Pi_{n,1}\big(C_{n,1}(I\varepsilon_{n,1})\big) = o(e^{-3n\varepsilon_{n,2}^2}).$$

## 4. Testing a Finite- versus an Infinite-dimensional Model

Suppose that there are two models, with the bigger models $\mathcal{P}_{n,1}$ infinite dimensional, and the alternative model a fixed parametric model $\mathcal{P}_{n,2} = \mathcal{P}_2 = \{p_\theta : \theta \in \Theta\}$, for $\Theta \subset \mathbb{R}^d$, equipped with a fixed prior $\Pi_{n,2} = \Pi_2$. Assume that $\lambda_{n,1} = \lambda_{n,2}$. We shall show that the Bayes factors are typically consistent in this situation: $\mathrm{BF}_n \to \infty$ if $p_0 \in \mathcal{P}_2$, and $\mathrm{BF}_n \to 0$ if $p_0 \notin \mathcal{P}_2$.

If the prior $\Pi_2$ is smooth in the parameter and the parametrization $\theta \mapsto p_\theta$ is regular, then, for any $\theta_0 \in \Theta$ and $\varepsilon \to 0$,

$$\Pi_2\Big(\theta : P_{\theta_0} \log \frac{p_{\theta_0}}{p_\theta} \leq \varepsilon^2, P_{\theta_0}\Big(\log \frac{p_{\theta_0}}{p_\theta}\Big)^2 \leq \varepsilon^2\Big) \sim C_{\theta_0} \varepsilon^d.$$

Therefore, if the true density $p_0$ is contained in $\mathcal{P}_2$, then $\Pi_{n,2}\big(B_{n,2}(\varepsilon_{n,2})\big) \gtrsim \varepsilon_{n,2}^d$, which is greater than $\exp[-n\varepsilon_{n,2}^2]$ for $\varepsilon_{n,2}^2 = D \log n/n$, for $D \geq d/2$ and sufficiently large $n$. (The logarithmic factor enters, because we use the crude prior mass condition (2.5) instead of the comparisons of prior mass in the main theorems, but it does not matter for this example.)

For this choice of $\varepsilon_{n,2}$ we have $\exp[n\varepsilon_{n,2}^2] = n^D$. Therefore, it follows from (2) of Corollary 3.1, that if $p_0 \in \mathcal{P}_2$, then the Bayes factor $\mathrm{BF}_n$ tends to $\infty$ as $n \to \infty$ as soon as there exists $\varepsilon_{n,1} > \varepsilon_{n,2}$ such that

$$\Pi_{n,1}\big(p : d(p, p_0) \leq I\varepsilon_{n1}\big) = o(n^{-3D}). \tag{4.1}$$

For an infinite-dimensional model $\mathcal{P}_{n,1}$ this is typically true, even if the models are nested, when $p_0$ is also contained in $\mathcal{P}_{n,1}$. In fact, we typically have for $p_0 \in \mathcal{P}_{n,1}$ that the left side is of the order $\exp[-Fn\varepsilon_{n,1}^2]$ for $\varepsilon_{n,1}$ the rate attached to the model $\mathcal{P}_{n,1}$. As for a true infinite-dimensional model this rate is not faster than $n^{-a}$ for some $a < 1/2$, this gives an upper bound of the order $\exp(-Fn^{1-2a})$, which is easily $o(n^{-3D})$. For $p_0$ not contained in the model $\mathcal{P}_{n,1}$, the prior mass in the preceding display will be even smaller than this.

If $p_0$ is not contained in the parametric model, then typically $d(p_0, \mathcal{P}_2) > 0$ and hence $d(p_0, \mathcal{P}_2) > I_n \varepsilon_{n,1}$ for any $\varepsilon_{n,1} \to 0$ and sufficiently slowly increasing $I_n$, as required in (1) of Corollary 3.1. To ensure that $\mathrm{BF}_n \to 0$, it suffices that for some $\varepsilon_{n,1} > \varepsilon_{n,2}$,

$$\Pi_{n,1}\Big(p: P_0 \log \frac{p_0}{p} \le \varepsilon_{n,1}^2, P_0\Big(\log \frac{p_0}{p}\Big)^2 \le \varepsilon_{n,1}^2\Big) \ge e^{-n\varepsilon_{n,1}^2}. \qquad (4.2)$$

This is the usual prior mass condition (cf. Ghosal et al. [2000]) for obtaining the rate of convergence $\varepsilon_{n,1}$ using the prior $\Pi_{n,1}$ on the model $\mathcal{P}_{n,1}$.

We present three concrete examples where the preceding can be made precise.

**Example 4.1** (Bernstein-Dirichlet mixtures). Bernstein polynomial densities with a Dirichlet prior on the mixing distribution and a geometric or Poisson prior on the order are described in Petrone [1999], Petrone and Wasserman [2002] and Ghosal [2001].

If we take these as the model $\mathcal{P}_{n,1}$ and prior $\Pi_{n,1}$, then the rate $\varepsilon_{n,1}$ is equal to $\varepsilon_{n,1} = n^{-1/3}(\log n)^{1/3}$ and (4.2) is satisfied, as shown in Ghosal [2001]. As the prior spreads its mass over an infinite-dimensional set that can approximate any smooth function, condition (4.1) will be satisfied for most true densities $p_0$. In particular, if $k_n$ is the minimal degree of a polynomial that is within Hellinger distance $n^{-1/3} \log n$ of $p_0$, then the left side of (4.1) is bounded by the prior mass of all Bernstein-Dirichlet polynomials of degree at least $k_n$, which is $e^{-ck_n}$ for some constant $c$ by construction. Thus (4.1) is certainly satisfied if $k_n \gg \log n$. Consequently, the Bayes factor is consistent for true densities that are not well approximable by polynomials.

**Example 4.2** (Log spline densities). Let $\mathcal{P}_{n,1}$ be equal to the set of log spline densities described in Section 5 of dimension $J \sim n^{1/(2\alpha+1)}$, equipped with the prior obtained by putting the uniform distribution on $[-M, M]^J$ on the coefficients. The corresponding rate can then be taken $\varepsilon_{n,1} = n^{-\alpha/(2\alpha+1)}\sqrt{\log n}$ (see Section 5.1). Conditions (4.1) and (4.2) can be verified easily by computations on the uniform prior, after translating the distances on the spline densities into the Euclidean distance on the coefficients (see Lemmas 7.4 and 7.6).

**Example 4.3** (Infinite dimensional normal model). Let $\mathcal{P}_{n,1}$ be the set of $N_\infty(\theta, I)$-distributions with $\theta = (\theta_1, \theta_2, \ldots)$ satisfying $\sum_{i=1}^\infty i^{2\alpha}\theta_i^2 < \infty$. (Thus a typical observation is an infinite sequence of independent normal variables with means $\theta_i$ and variances 1.) Equip it with the prior obtained by letting the $\theta_i$ be independent Gaussian variables with mean 0 and variances $i^{-(2q+1)}$. Take $\mathcal{P}_{n,2}$ equal to the submodel indexed by all $\theta$ with $\theta_i = 0$ for all $i \ge 2$, equipped with a positive smooth (for instance Gaussian) prior on $\theta_1$. This model is equivalent to the signal plus Gaussian white noise model, for which Bayesian procedures were studied in Freedman [1999], Zhao [2000], Belitser and Ghosal [2003] and Ghosal and van der Vaart [2007a].

The Kullback-Leibler and squared Hellinger distances on $\mathcal{P}_{n,1}$ are, up to constants, essentially equivalent to the squared $\ell_2$-distance on the parameter $\theta$, when the distance is bounded (see e.g. Lemma 6.1 of Belitser and Ghosal [2003]). This allows to verify (4.1) by calculations on Gaussian variables, after truncating the parameter set. For a sufficiently large constant $M$, consider sieves $\mathcal{P}'_{n,1} = \{\theta: \sum_{i=1}^{\infty} i^{2q'}\theta_i^2 \leq M\}$ for some $q' < q$, and $\mathcal{P}'_{n,2} = \{|\theta_1| \leq M\}$, respectively. The posterior probabilities of the complements of these sieves are small in probability respectively by Lemma 3.2 of Belitser and Ghosal [2003] and Lemma 7.2 of Ghosal et al. [2000]. Hence in view of Remark 2.1 it suffices to perform the calculations on $\mathcal{P}'_{n,1}$ and $\mathcal{P}'_{n,2}$.

By Lemmas 6.3 and 6.4 of Belitser and Ghosal [2003] it follows that the conditions of (2) of Corollary 3.1 holds for $\varepsilon_{n,1} = \max(n^{-q'/(2q'+1)}, n^{-q/(2q+1)}) = n^{-q'/(2q'+1)}$, provided that (4.1) can be verified. Now, for any $\theta_0 \in \ell_2$,

$$\Pi_{n,1}\big(\theta: \|\theta - \theta_0\| \leq \varepsilon_{n,1}\big) \leq \prod_{i=1}^{\infty} \Pi_{n,1}\big(|\theta_i - \theta_{i0}| \leq \varepsilon_{n,1}\big) \leq \prod_{i=1}^{\infty} \big(2\Phi(i^{q+1/2}\varepsilon_{n,1}) - 1\big).$$

For $i \leq n^{2q'/((2q+1)(2q'+1))}$ the argument in the normal distribution function is bounded above by 1, and then the corresponding factor is bounded by $2\Phi(1) - 1 < 1$. It follows that the right side of the last display is bounded by a term of the order $e^{-c'n^{2q'/((2q+1)(2q'+1))}}$, for some positive constant $c$. This easily shows that (4.1) is satisfied.

## 5. Log Spline Models

Log spline density models, introduced in Stone [1990], are exponential families constructed as follows.

For a given "resolution" $K \in \mathbb{N}$ partition the half open unit interval $[0, 1)$ into $K$ subintervals $\big[(k-1)/K, k/K\big)$ for $k = 1, \ldots, K$. The linear space of *splines* of "order" $q \in \mathbb{N}$ relative to this partition is the set of all continuous functions $f: [0, 1] \to \mathbb{R}$ that are $q - 2$ times differentiable on $[0, 1)$ and whose restriction to every of the partitioning intervals $\big[(k-1)/K, k/K\big)$ is a polynomial of degree strictly less than $q$. It can be shown that these splines form a $J = q + K - 1$-dimensional vector space. A convenient basis is the set of B-splines $B_{J,1}, \ldots, B_{J,J}$, defined e.g. in de Boor [2001]. The exact nature of these functions does not matter to us here, but the following properties are essential (cf. de Boor [2001], pp 109–110):

- $B_{J,j} \geq 0, \qquad j = 1, \ldots, J$
- $\sum_{j=1}^{J} B_{J,j} \equiv 1$
- $B_{J,j}$ is supported on an interval of length $q/K$
- at most $q$ functions $B_{J,j}$ are nonzero at every given $x$.

The first two properties express that the basis elements form a partition of unity, and the third and fourth property mean that their supports are close to being disjoint if $K$ is very large relative to $q$. This renders the B-spline basis stable for numerical computation, and also explains the simple inequalities between norms of linear combinations of the basis functions and norms of the coefficients given in Lemma 7.2 below.

For $\theta \in \mathbb{R}^J$ let $\theta^T B_J = \sum_j \theta_j B_{J,j}$ and define

$$p_{J,\theta}(x) = e^{\theta^T B_J(x) - c_J(\theta)}, \qquad e^{c_J(\theta)} = \int_0^1 e^{\theta^T B_J(x)} \, dx.$$

Thus $p_{J,\theta}$ is a probability density that belongs to a $J$-dimensional exponential family with sufficient statistics the B-spline functions. Since the B-splines add up to unity, the family is actually of dimension $J - 1$ and we can restrict $\theta$ to the subset of $\theta \in \mathbb{R}^J$ such that $\theta^T 1 = 0$.

Splines possess excellent approximation properties for smooth functions, where the error is smaller if the function is smoother or the dimension of the spline space is higher. More precisely, a function $f \in C^\alpha[0,1]$ can be approximated with an error of order $(1/J)^\alpha$ by splines of order $q \geq \alpha$ and dimension $J$. Because there are $J - 1$ free base coefficients, the variance of a best estimate in a $J$-dimensional spline space can be expected to be of order $J/n$. Therefore, we may expect to determine an optimal dimension $J$ for a given smoothness level $\alpha$ from the bias-variance trade-off $J/n \sim (1/J)^{2\alpha}$. This leads to the dimension $J_{n,\alpha} \sim n^{1/(2\alpha+1)}$, and the "usual" rate of convergence $n^{-\alpha/(2\alpha+1)}$.

This informal calculation was justified for maximum likelihood and Bayesian estimators in Stone [1990] and Ghosal et al. [2000], respectively. The paper Stone [1990] showed that the maximum likelihood estimator of $p$ in the model $\{p_{J,\theta} : J = J_{n,\alpha}, \theta \in \mathbb{R}^{J_{n,\alpha}}, \theta^T 1 = 0\}$ achieves the rate of convergence $n^{-\alpha/(2\alpha+1)}$ if the true density $p_0$ belongs to $C^\alpha[0,1]$. The paper Ghosal et al. [2000] showed that a Bayes procedure with deterministic dimension $J_{n,\alpha}$ and a smooth prior on the coefficients $\theta \in \mathbb{R}^{J_{n,\alpha}}$ achieves the same (posterior) rate. (In both papers it is assumed that the true density is also bounded away from zero.)

Both the maximum likelihood estimator and the Bayesian estimator described previously depend on $\alpha$. They can be made rate-adaptive to $\alpha$ by a variety of means. We shall consider several Bayesian schemes, based on different choices of priors $\bar{\Pi}_{n,\alpha}$ on the coefficients and $\lambda_n$ on the dimensions $J_{n,\alpha}$ of the spline spaces. Thus $\bar{\Pi}_{n,\alpha}$ will be a prior on $\mathbb{R}^{J_{n,\alpha}}$ for

$$J_{n,\alpha} = \lfloor n^{1/(2\alpha+1)} \rfloor, \tag{5.1}$$

the prior $\Pi_{n,\alpha}$ on densities will be the distribution induced under the map $\theta \mapsto p_{J,\theta}$, where $J = J_{n,\alpha}$, and $\lambda_n$ is a prior on the regularity parameter $\alpha$.

We always choose the order of the splines involved in the construction of the $\alpha$th log spline model at least $\alpha$.

We shall assume that the true density $p_0$ is bounded away from zero next to being smooth, so that the Hellinger and $L_2$-metrics are equivalent. We shall in fact assume that uniform upper and lower bounds are known, and construct the priors on sets of densities that are bounded away from zero and infinity. It follows from Lemma 7.3 that the latter is equivalent to restricting the coefficient vector $\theta$ in $p_{J,\theta}$ to a rectangle $[-M, M]^J$ for some $M$. We shall construct our priors on this rectangle, and assume that the true density $p_0$ is within the range of the corresponding spline densities, i.e. $\|\log p_0\|_\infty \leq \underline{C}_4 M$ for the constant $\underline{C}_4$ of Lemma 7.3. Extension to unknown $M$ through a second shell of adaptation is possible (see e.g. Lember and van der Vaart [2007]), but will not be pursued here.

In the next three sections we discuss three examples of priors. In the first example we combine *smooth* priors $\bar{\Pi}_{n,\alpha}$ on the coefficients with fixed model weights $\lambda_{n,\alpha} = \lambda_\alpha$. These natural priors lead to adaptation up to a logarithmic factor. Even though we only prove an upper bound, we believe that the logarithmic factor is not a defect of our proof, but connected to this prior. In the second example we show how the logarithmic factor can be removed by using special model weights $\lambda_{n,\alpha}$ that put *less* mass on small models. In the third example we show that the logarithmic factor can also be removed by using *discrete* priors $\Pi_{n,\alpha}$ for the coefficients combined with model weights that put *more* mass on small models. One might conclude from these examples that the fine details of rates depend on a delicate balance between model priors and model weights. The good news is that all three priors work reasonably well.

### 5.1. Flat priors

In this section we consider prior distributions $\bar{\Pi}_{n,\alpha}$ that possess Lebesgue densities on $\{\theta \in \mathbb{R}^{J_{n,\alpha}} : \theta^T 1 = 0\}$ that vanish outside a big block $[-M, M]^{J_{n,\alpha}}$ and are bounded above and below by $d^{J_{n,\alpha}}$ and $D^{J_{n,\alpha}}$, respectively, for given constants $0 < d \leq D < \infty$. We combine these with fixed prior weights $\lambda_{n,\alpha} = \mu_\alpha > 0$ on the regularity parameter, which we restrict to $A = \{\alpha \in \mathbb{Q}^+ : \alpha \geq \underline{\alpha}\}$ for some (known) constant $\underline{\alpha} > 0$ and assume to satisfy $\sum_{\alpha \in A} \sqrt{\mu_\alpha} < \infty$.

**Theorem 5.1.** *If $p_0 \in C^\beta[0, 1]$ for some $\beta \in \mathbb{Q}^+ \cap [\underline{\alpha}, \infty)$ and $\|\log p_0\|_\infty < \underline{C}_4 M$, then there exist a constant $B$ such that $P_0^n \Pi_n \big( p : \|p - p_0\|_2 \geq B\varepsilon_{n,\beta} \big| X_1, \ldots, X_n \big) \to 0$, for $\varepsilon_{n,\beta} = n^{-\beta/(2\beta+1)} \sqrt{\log n}$.*

*Proof.* The dimension numbers $J_{n,\alpha}$ defined in (5.1) relate to the present rates $\varepsilon_{n,\alpha} := n^{-\alpha/(2\alpha+1)} \sqrt{\log n}$ as $J_{n,\alpha} \log n \sim n\varepsilon_{n,\alpha}^2$.

By Lemma 7.7 condition (2.1) is satisfied for any $\varepsilon_{n,\alpha}$ such that $n\varepsilon_{n,\alpha}^2 \gtrsim J_{n,\alpha}$, and hence certainly for the present $\varepsilon_{n,\alpha}$. The constants $E_\alpha$ do not depend on $\alpha$, and hence both $E$ and $\underline{E}$ in Theorem 2.1 can be taken equal to a single constant $E$.

Because $\|\log p_0\|_\infty < \underline{C}_4 M$ by assumption, the Hellinger distance of $p_0$ to $\mathcal{P}_J$ is bounded above by a multiple of $J^{-\beta}$ by Lemma 7.8. By Lemma 7.6 for $\varepsilon_{n,\beta} \gtrsim J_{n,\beta}^{-\beta}$, some constants $\underline{A}$ and $\overline{A}$, and sufficiently large $n$,

$$\Pi_{n,\alpha}\big(C_{n,\alpha}(\varepsilon)\big) \leq \overline{\Pi}_{n,\alpha}\big(\theta \in \Theta_{J_{n,\alpha}} \colon \|\theta - \theta_{J_{n,\alpha}}\|_2 \leq \underline{A}\sqrt{J_{n,\alpha}}\,\varepsilon\big),$$

$$\Pi_{n,\beta}\big(B_{n,\beta}(\varepsilon_{n,\beta})\big) \geq \overline{\Pi}_{n,\beta}\big(\theta \in \Theta_{J_{n,\beta}} \colon \|\theta - \theta_{J_{n,\beta}}\|_2 \leq 2\overline{A}\sqrt{J_{n,\beta}}\,\varepsilon_{n,\beta}\big).$$

Because $\theta_J \in \Theta_J$ by its definition, the set $\{\theta \in \Theta_J \colon \|\theta - \theta_J\|_2 \leq \varepsilon\}$ contains at least a fraction $2^{-J}$ of the volume of the ball of radius $\varepsilon$ around $\theta_J$, even though it does not contain the full ball if $\theta_J$ is near the boundary of $\Theta_J$. It follows that, for any $\alpha, \beta \in A$ and $\varepsilon$, and $v_J$ the volume of the $J$-dimensional Euclidean ball,

$$\frac{\Pi_{n,\alpha}\big(C_{n,\alpha}(i\varepsilon)\big)}{\Pi_{n,\beta}\big(B_{n,\beta}(\varepsilon_{n,\beta})\big)} \leq \frac{(\underline{A}Di\varepsilon\sqrt{J_{n,\alpha}})^{J_{n,\alpha}} v_{J_{n,\alpha}}}{(\overline{A}d\varepsilon_{n,\beta}\sqrt{J_{n,\beta}})^{J_{n,\beta}} v_{J_{n,\beta}}} \lesssim \frac{(\overline{a}i\varepsilon)^{J_{n,\alpha}}}{(a\varepsilon_{n,\beta})^{J_{n,\beta}}}, \qquad (5.2)$$

for suitable constants $a$ and $\overline{a}$, in view of Lemma 7.9.

If $\alpha < \beta$, then with $\varepsilon = \varepsilon_{n,\alpha}$, inequality (5.2) yields

$$\frac{\Pi_{n,\alpha}\big(C_{n,\alpha}(i\varepsilon_{n,\alpha})\big)}{\Pi_{n,\beta}\big(B_{n,\beta}(\varepsilon_{n,\beta})\big)} \lesssim \frac{(\overline{a}i\varepsilon_{n,\alpha})^{J_{n,\alpha}}}{(a\varepsilon_{n,\beta})^{J_{n,\beta}}}$$

$$= \exp\Big[J_{n,\alpha}\Big(\log(\overline{a}i\varepsilon_{n,\alpha}) - \frac{J_{n,\beta}}{J_{n,\alpha}}\log(a\varepsilon_{n,\beta})\Big)\Big]$$

$$\leq \exp\Big[J_{n,\alpha}\Big(\log(\overline{a}i) + \frac{1}{H}|\log a| + \log\varepsilon_{n,\alpha} - \frac{1}{H}\log\varepsilon_{n,\beta}\Big)\Big],$$

because $J_{n,\alpha} > HJ_{n,\beta}$ for $\alpha < \beta$. Here,

$$\log\varepsilon_{n,\alpha} - \frac{1}{H}\log\varepsilon_{n,\beta} = \Big(\frac{1}{H}\frac{\beta}{2\beta+1} - \frac{\alpha}{2\alpha+1}\Big)\log n + \frac{1}{2}\Big(1 - \frac{1}{H}\Big)\log\log n.$$

For sufficiently large $H$ the coefficient of $\log n$ is negative, uniformly in $\alpha \geq \underline{\alpha}$. Condition (2.2) is easily satisfied for such $H$, with $\mu_{n,\alpha} = \mu_\alpha/\mu_\beta$ and arbitrarily small $L > 0$.

If $\alpha < \beta$, then with $\varepsilon = IB\varepsilon_{n,\alpha}$, inequality (5.2) and similar calculations yield

$$e^{2n\varepsilon_{n,\beta}^2}\frac{\Pi_{n,\alpha}\big(C_{n,\alpha}(IB\varepsilon_{n,\alpha})\big)}{\Pi_{n,\beta}\big(B_{n,\beta}(\varepsilon_{n,\beta})\big)}$$

$$\lesssim \exp\Big[J_{n,\alpha}\Big(\log(\overline{a}IB) + \log\varepsilon_{n,\alpha} - \frac{J_{n,\beta}}{J_{n,\alpha}}\log(a\varepsilon_{n,\beta}) + 2\frac{J_{n,\beta}}{J_{n,\alpha}}\log n\Big)\Big]$$

$$\leq \exp\Big[J_{n,\alpha}\Big(\log(\overline{a}IB) + \frac{1}{H}|\log a| + \log\varepsilon_{n,\alpha} - \frac{1}{H}\log\varepsilon_{n,\beta} + 2\frac{1}{H}\log n\Big)\Big].$$

By the same arguments as before for sufficiently large $H$ the exponent is smaller than $-J_{n,\alpha}c\log n$ for a positive constant $c$, uniformly in $\alpha \geq \underline{\alpha}$, eventually. This implies that (2.4) is fulfilled.

With $\varepsilon = \varepsilon_{n,\beta}$, inequality (5.2) yields

$$\frac{\Pi_{n,\alpha}\big(C_{n,\alpha}(i\varepsilon_{n,\beta})\big)}{\Pi_{n,\beta}\big(B_{n,\beta}(\varepsilon_{n,\beta})\big)} \lesssim \frac{(\overline{a}i\varepsilon_{n,\beta})^{J_{n,\alpha}}}{(a\varepsilon_{n,\beta})^{J_{n,\beta}}}$$

$$= \exp\Big[J_{n,\beta}\Big(\frac{J_{n,\alpha}}{J_{n,\beta}}\big(\log(\overline{a}i) + \log\varepsilon_{n,\beta}\big) - \log(a\varepsilon_{n,\beta})\Big)\Big].$$

If $\alpha \gtrsim \beta$ the right side is bounded above by

$$\exp\big[J_{n,\beta}\big(H|\log(\overline{a}i)| - \log(a\varepsilon_{n,\beta})\big)\big] \le \exp\big[J_{n,\beta}(\log n)Li^2\big],$$

for sufficiently large $i$, where $L$ may be an arbitrarily small constant. Hence condition (2.3) is fulfilled.

The theorem is a consequence of Theorem 2.1, with $A_n$ of this theorem equal to the present $A$. ∎

### 5.2. Flat priors, decreasing weights

The constant weights $\lambda_{n,\alpha} = \mu_\alpha$ used in the preceding subsection resulted in an additional logarithmic factor in the rate. The following theorem shows that for $A = \{\alpha_1, \alpha_2, \ldots, \alpha_N\}$ a finite set, this factor can be removed by choosing the weights

$$\lambda_{n,\alpha} \propto \prod_{\gamma \in A : \gamma < \alpha} (C\varepsilon_{n,\gamma})^{J_{n,\gamma}}. \tag{5.3}$$

These weights are decreasing in $\alpha$, unlike the weights in (2.8). Thus the present prior puts less weight on the smaller, more regular models. We use the same priors $\Pi_{n,\alpha}$ on the spline models as in Section 5.1.

**Theorem 5.2.** *Let $A$ be a finite set. If $p_0 \in C^\beta[0,1]$ for some $\beta \in A$ and $\|\log p_0\|_\infty < \underline{C}_4 M$ and sufficiently large $C$, then there exist a constant $B$ such that $P_0^n \Pi_n\big(p : \|p - p_0\|_2 \ge B\varepsilon_{n,\beta}|X_1, \ldots, X_n\big) \to 0$ for $\varepsilon_{n,\beta} = n^{-\beta/(2\beta+1)}$.*

*Proof.* Let $\varepsilon_{n,\alpha} = n^{-\alpha/(2\alpha+1)}$, so that $J_{n,\alpha} \sim n\varepsilon_{n,\alpha}^2$, and $J_{n,\alpha'}/J_{n,\alpha} \ll n^{-c}$ for some $c > 0$ whenever $\alpha' > \alpha$. Assume without loss of generality that $A = \{\alpha_1, \ldots, \alpha_N\}$ is indexed in its natural order.

If $r < s$, and hence $\alpha_r < \alpha_s$, then, by inequality (5.2),

$$\frac{\lambda_{n,\alpha_r}\Pi_{n,\alpha_r}\big(C_{n,\alpha_r}(i\varepsilon_{n,\alpha_r})\big)}{\lambda_{n,\alpha_s}\Pi_{n,\alpha_s}\big(B_{n,\alpha_s}(\varepsilon_{n,\alpha_s})\big)}$$

$$\lesssim \exp\Big[J_{n,\alpha_r}\Big(\log(\overline{a}i\varepsilon_{n,\alpha_r}) - \frac{J_{n,\alpha_s}}{J_{n,\alpha_r}}\log(a\varepsilon_{n,\alpha_s}) - \sum_{k=r}^{s-1}\frac{J_{n,\alpha_k}}{J_{n,\alpha_r}}\log(C\varepsilon_{n,\alpha_k})\Big)\Big]$$

$$= \exp\Big[J_{n,\alpha_r}\Big(\log\Big(\frac{\overline{a}i}{C}\Big) - \frac{J_{n,\alpha_s}}{J_{n,\alpha_r}}\log(a\varepsilon_{n,\alpha_s}) - \sum_{k=r+1}^{s-1}\frac{J_{n,\alpha_k}}{J_{n,\alpha_r}}\log(C\varepsilon_{n,\alpha_k})\Big)\Big].$$

The exponent takes the form $J_{n,\alpha_r}\big(\log(\overline{a}i/C) + o(1)\big)$. Applying this with $\alpha_r = \alpha < \beta = \alpha_s$, we conclude that (2.2) holds for every $C$ and $\mu_{n,\alpha} = 1$, for sufficiently large $i$, for an arbitrarily small constant $L$, eventually.

Similarly, again with $\alpha_r = \alpha < \beta = \alpha_s$,

$$e^{2n\varepsilon_{n,\alpha_s}^2} \frac{\lambda_{n,\alpha_r}\Pi_{n,\alpha_r}\big(C_{n,\alpha_r}(IB\varepsilon_{n,\alpha_r})\big)}{\lambda_{n,\alpha_s}\Pi_{n,\alpha_s}\big(B_{n,\alpha_s}(\varepsilon_{n,\alpha_s})\big)}$$

$$= \exp\Big[J_{n,\alpha_r}\Big(\log\Big(\frac{\overline{a}IB}{C}\Big) - \frac{J_{n,\alpha_s}}{J_{n,\alpha_r}}\log(a\varepsilon_{n,\alpha_s})$$

$$- \sum_{k=r+1}^{s-1} \frac{J_{n,\alpha_k}}{J_{n,\alpha_r}}\log(C\varepsilon_{n,\alpha_k}) + \frac{2J_{n,\alpha_s}}{J_{n,\alpha_r}}\Big)\Big].$$

This tends to 0 if $C > \overline{a}IB$. Hence, for $C$ big enough, condition (2.4) is fulfilled as well.

Finally, choose $\alpha_r = \beta < \alpha = \alpha_s$ and note that

$$\frac{\lambda_{n,\alpha_s}\Pi_{n,\alpha_s}\big(C_{n,\alpha_s}(i\varepsilon_{n,\alpha_r})\big)}{\lambda_{n,\alpha_r}\Pi_{n,\alpha_r}\big(B_{n,\alpha_r}(\varepsilon_{n,\alpha_r})\big)}$$

$$\lesssim \exp\Big[J_{n,\alpha_r}\Big(\frac{J_{n,\alpha_s}}{J_{n,\alpha_r}}\big(\log(\overline{a}i\varepsilon_{n,\alpha_r}) - \log(a\varepsilon_{n,\alpha_r}) + \sum_{k=r}^{s-1}\frac{J_{n,\alpha_k}}{J_{n,\alpha_r}}\log(C\varepsilon_{n,\alpha_k})\big)\Big)\Big]$$

$$= \exp\Big[J_{n,\alpha_r}\Big(\frac{J_{n,\alpha_s}}{J_{n,\alpha_r}}\log(\overline{a}i\varepsilon_{n,\alpha_r}) + \log\Big(\frac{C}{a}\Big) + \sum_{k=r+1}^{s-1}\frac{J_{n,\alpha_k}}{J_{n,\alpha_r}}\log(C\varepsilon_{n,\alpha_k})\Big)\Big].$$

Here the exponent is of the order $J_{n,\alpha_r}\big(\log(C/a) + o(1)\log i + o(1)\big)$. We conclude that the condition (2.3) holds.

The theorem is a consequence of Theorem 2.1, with $A_n$ of this theorem equal to the present $A$. ∎

The proof of Theorem 5.2 relies on the fact that $J_{n,\alpha_r} \ll J_{n,\alpha_s}$ if $s < r$, and for that reason it does not extend to the more general case of a dense set $A$ as considered in Theorem 5.1. On the other hand, it would be possible to extend the theorem to a countable totally ordered set $\alpha_1 < \alpha_2 < \cdots$ by using the weights (5.3) restricted to sets $\alpha_1 < \alpha_2 < \cdots < \alpha_{M_n}$ for $M_n \uparrow \infty$.

The existence of rate-adaptive priors that yield the rate without log-factor in the general countable case is subject for further research. There are some reasons to believe that this task is achievable with some more elaborate priors as these in (5.3). For example, one could consider more general priors than (5.3) of the type

$$\lambda_{n,\alpha} \propto \prod_{\gamma \in A_n} (C_{\gamma,\alpha}\varepsilon_{n,\gamma})^{\lambda_{\gamma,\alpha}J_{n,\gamma}}.$$

The truncation-set $A_n$ as well as the constants $C_{\gamma,\alpha}$ and $\lambda_{\gamma,\alpha}$ must be carefully chosen.

### *5.3. Discrete priors, increasing weights*

In this section we choose the priors $\Pi_{n,\alpha}$ to be discrete on a suitable subset of the $J_{n,\alpha}$-dimensional log spline densities, constructed as follows.

According to Kolmogorov and Tihomirov [1961] (cf. Theorem 2.7.1 in van der Vaart and Wellner [1996]) the unit ball $C_1^\alpha[0,1]$ of the Hölder space $C^\alpha[0,1]$ has entropy $\log N\left(\varepsilon, C_1^\alpha[0,1], \|\cdot\|_\infty\right)$ of the order $(1/\varepsilon)^{1/\alpha}$, as $\varepsilon \downarrow 0$, relative to the uniform norm. Then it follows that there exists a set of $N_{n,\alpha} \lesssim (M/\varepsilon_{n,\alpha})^{1/\alpha}$ functions $f_1, \ldots, f_{N_{n,\alpha}}$ such that every $f$ with Hölder norm smaller than a given constant $M$ is within uniform distance $\varepsilon_{n,\alpha} := n^{-\alpha/(2\alpha+1)}$ of some $f_i$. These functions can without loss of generality be chosen with Hölder norm bounded by $M$. By the approximation properties of spline spaces (cf. Lemma 7.1), we can find $\theta_i \in \mathbb{R}^{J_{n,\alpha}}$ such that $\|\theta_i^T B_{J_{n,\alpha}} - f_i\|_\infty \le C_{q,\alpha}(1/J_{n,\alpha})^\alpha$. Define $\bar{\Pi}_{n,\alpha}$ to be the uniform probability distribution on the collection $\theta_1, \ldots, \theta_{N_{n,\alpha}}$.

We combine the resulting prior $\Pi_{n,\alpha}$ on log spline densities with model weights on the index set $A = \mathbb{Q}^+$ of the form (2.8), where $(\mu_\alpha : \alpha \in \mathbb{Q}^+)$ is a strictly positive measure with $\sum_{\alpha \in \mathbb{Q}^+} \sqrt{\mu_\alpha} < \infty$ and $C$ is an arbitrary positive constant.

**Theorem 5.3.** *If $p_0 \in C^\beta[0,1]$ for some $\beta \in A$ and $\|\log p_0\|_\beta < M$, then there exist a constant $B$ such that $P_0^n \Pi_n\left(p : h(p, p_0) \ge B\varepsilon_{n,\beta} | X_1, \ldots, X_n\right) \to 0$ for $\varepsilon_{n,\beta} = n^{-\beta/(2\beta+1)}$.*

*Proof.* By construction there exists an element $\theta_i$ in the support of $\bar{\Pi}_{n,\beta}$ such that

$$\|\log p_0 - \theta_i^T B_{J_{n,\beta}}\|_\infty \lesssim \varepsilon_{n,\beta} + (1/J_{n,\beta})^\beta \lesssim \varepsilon_{n,\beta}.$$

It follows that the function $e^{\theta_i^T B_{J_{n,\beta}}}$ is sandwiched between the functions $p_0 e^{-d\varepsilon_{n,\beta}}$ and $p_0 e^{+d\varepsilon_{n,\beta}}$, for some constant $d$. Consequently, the norming constant satisfies $|e^{-c(\theta_i)} - 1| \lesssim \varepsilon_{n,\beta}$, and hence $\|p_0 - p_{J_{n,\beta},\theta_i}\|_\infty \lesssim \varepsilon_{n,\beta}$. Because $p_0$ is bounded away from zero and infinity, this implies that $p_{J_{n,\beta},\theta_i}$ is in the Kullback-Leibler neighbourhood $B_{n,\beta}(D\varepsilon_{n,\beta})$ of $p_0$, for some constant $D$ (cf. Lemma 8 in Ghosal and van der Vaart [2007b]). Because $\Pi_{n,\beta}$ is the uniform measure on the $N_{n,\beta}$ log spline densities of this type, it follows that $\Pi_{n,\beta}\left(B_{n,\beta}(D\varepsilon_{n,\beta})\right) \ge N_{n,\beta}^{-1} \ge \exp[-Fn\varepsilon_{n,\beta}^2]$, for some positive constant $F$.

In view of (2.8) it follows, for any $\varepsilon$,

$$\frac{\lambda_{n,\alpha}}{\lambda_{n,\beta}} \frac{\Pi_{n,\alpha}\left(C_{n,\alpha}(i\varepsilon)\right)}{\Pi_{n,\beta}\left(B_{n,\beta}(D\varepsilon_{n,\beta})\right)} \le \frac{\lambda_{n,\alpha}}{\lambda_{n,\beta}} e^{FJ_{n,\beta}} = \frac{\mu_\alpha}{\mu_\beta} e^{-CJ_{n,\alpha}} e^{(F+C)J_{n,\beta}}.$$

Define the sets of indices $\alpha < \beta$ and $\alpha \gtrsim \beta$ as in Theorem 2.1, relative to a given constant $H$. Thus $\alpha < \beta$ is equivalent to $J_{n,\alpha} > HJ_{n,\beta}$ and hence the sum over $\alpha < \beta$ of the preceding display can be bounded above by

$$e^{(C-CH/2+F)J_{n,\beta}} \sum_{\alpha < \beta} \frac{\mu_\alpha}{\mu_\beta} e^{-CJ_{n,\alpha}/2}.$$

The leading term is $o(e^{-2n\varepsilon_{n,\beta}^2})$ provided $H$ is big enough, and the sum is bounded by assumption. Thus (2.4) is fulfilled for any constant $B$. Furthermore, condition (2.2) holds trivially with $\mu_{n,\alpha} = (\mu_\alpha/\mu_\beta)e^{-CJ_{n,\alpha}/2}$. Condition (2.3) clearly holds for sufficiently large $i$, with the same choice of $\mu_{n,\alpha}$ and any $L > 0$.

The theorem is a consequence of Theorem 2.1, with $A_n$ of this theorem equal to the present $A$. ∎

## 6. Proof of the main theorems

We start by extending results from Ghosal and van der Vaart [2007b], Ghosal et al. [2000], LeCam [1973] and Birgé [1983] on the existence of tests of certain tests under local entropy of a statistical model. The results differ from the last three references by inclusion of weights $\alpha$ and $\beta$; relative to Ghosal and van der Vaart [2007b] the difference is the use of local rather than global entropy.

Let $d$ be a metric which induces convex balls and is bounded above on $\mathcal{P}$ by the Hellinger metric $h$.

**Lemma 6.1.** *For any dominated convex set of probability measures $\mathcal{P}$, and any constants $\alpha, \beta > 0$ and any $n$ there exists a test $\phi$ with*

$$\sup_{Q \in \mathcal{P}} \left( \alpha P_0^n \phi + \beta Q^n (1 - \phi) \right) \leq \sqrt{\alpha \beta} \, e^{-\frac{1}{2}nh^2(P_0, \mathcal{P})}.$$

*Proof.* This follows by minor adaptation of a result of Le Cam [1986]. The essence is that, by the minimax theorem,

$$\inf_\phi \sup_{Q \in \mathcal{P}} \left( \alpha P_0^n \phi + \beta Q^n (1 - \phi) \right)$$

$$= \sup_{Q_n \in \text{conv}\,(\mathcal{P}^n)} \inf_\phi \left( \alpha P_0^n \phi + \beta Q_n (1 - \phi) \right)$$

$$= \sup_{Q_n \in \text{conv}\,(\mathcal{P}^n)} \left( \alpha P_0^n (\alpha p_0^n < \beta q_n) + \beta Q_n (\alpha p_0^n > \beta q_n) \right)$$

$$\leq \sup_{Q_n \in \text{conv}\,(\mathcal{P}^n)} \left( \alpha \int_{\alpha p_0^n < \beta q_n} \sqrt{p_0^n} \sqrt{(\beta/\alpha)q_n} + \beta \int_{\alpha p_0^n > \beta q_n} \sqrt{(\alpha/\beta)p_0^n} \sqrt{q_n} \right)$$

$$\leq \sup_{Q_n \in \text{conv}\,(\mathcal{P}^n)} \sqrt{\alpha \beta} \int \sqrt{p_0^n} \sqrt{q_n}.$$

Next we use the convexity of $\mathcal{P}$ to see that this is bounded above by (see Le Cam [1986] or Lemma 6.2 in Kleijn and van der Vaart [2006])

$$\sqrt{\alpha \beta} \left( \sup_{Q \in \mathcal{P}} \int \sqrt{p_0} \sqrt{q} \right)^n.$$

Finally we express the affinity $\int \sqrt{p_0} \sqrt{q}$ in the Hellinger distance as $1 - \frac{1}{2}h^2(P_0, Q)$ and use the inequality $1 - x \leq e^{-x}$, for $x > 0$. ∎

**Corollary 6.1.** *For any dominated set of probability measures $\mathcal{P}$ with $d(P_0, \mathcal{P}) \geq 3\varepsilon$, any $\alpha, \beta > 0$ and any $n$ there exists a test $\phi$ with*

$$P_0^n \phi \leq \sqrt{\frac{\beta}{\alpha}} \, N(\varepsilon, \mathcal{P}, d) \, e^{-n\varepsilon^2},$$

$$\sup_{Q \in \mathcal{P}} Q^n (1 - \phi) \leq \sqrt{\frac{\alpha}{\beta}} \, e^{-n\varepsilon^2}.$$

*Proof.* Choose a maximal $2\varepsilon$-separated set $\mathcal{P}'$ of points in $\mathcal{P}$. Then the balls $B_{Q'}$ of radius $2\varepsilon$ centered at the points in $\mathcal{P}'$ cover $\mathcal{P}$, whence their number is bounded by $N(2\varepsilon, \mathcal{P}, d)$. Furthermore, these balls are convex by assumption, and are at distance $3\varepsilon - \varepsilon = 2\varepsilon$ from $P_0$. The latter is true both for the distance $d$ and the Hellinger distance, which is larger by assumption. For every ball $B_{Q'}$ attached to a point $Q' \in \mathcal{P}'$ there exists a test $\omega_{Q'}$ with the properties as in Lemma 6.1 with $\mathcal{P}$ taken equal to $B_{Q'}$. Let $\phi$ be the maximum of all tests attached in this way to some point $Q \in \mathcal{P}'$. Then

$$P_0^n \phi \leq \sum_{Q \in \mathcal{P}'} P_0^n \omega_{Q'} \leq \sum_{Q \in \mathcal{P}'} \sqrt{\frac{\beta}{\alpha}} \, e^{-n2\varepsilon^2},$$

$$\sup_{Q \in \mathcal{P}} Q^n (1 - \phi) \leq \sup_{Q \in \mathcal{P}'} Q^n (1 - \omega_{Q'}) \leq \sqrt{\frac{\alpha}{\beta}} \, e^{-n2\varepsilon^2}.$$

The right sides can be further bounded as desired. ∎

**Lemma 6.2.** *Suppose that for a dominated set of probability measures $\mathcal{P}$, some nonincreasing function $\varepsilon \mapsto N(\varepsilon)$, some $\varepsilon_0 \geq 0$, and for every $\varepsilon > \varepsilon_0$,*

$$N\Big(\frac{\varepsilon}{3}, \big\{ p \in \mathcal{P} \colon \varepsilon \leq d(p, p_0) \leq 2\varepsilon \big\}, d\Big) \leq N(\varepsilon).$$

*Then for every $\varepsilon > \varepsilon_0$ and every $\alpha, \beta > 0$ there exists a test $\phi$ (depending on $\varepsilon$ but not on $i$) such that for every $i \in \mathbb{N}$,*

$$P_0^n \phi \leq \sqrt{\frac{\beta}{\alpha}} N(\varepsilon) e^{-n\varepsilon^2/9} \frac{1}{1 - e^{-n\varepsilon^2/9}},$$

$$\sup_{p \in \mathcal{P} \colon d(p, p_0) > i\varepsilon} P^n (1 - \phi) \leq \sqrt{\frac{\alpha}{\beta}} e^{-n\varepsilon^2 i^2/9},$$

*Proof.* For $j \in \mathbb{N}$ let $\mathcal{P}_j = \{ p \in \mathcal{P} \colon j\varepsilon < d(p, p_0) \leq (j+1)\varepsilon \}$. Because the set $\mathcal{P}_j$ has distance $3(j\varepsilon/3)$ to $p_0$, the preceding corollary implies the existence of a test $\phi_j$ with

$$P_0^n \phi_j \leq \sqrt{\frac{\beta}{\alpha}} N(j\varepsilon) e^{-nj^2\varepsilon^2/9},$$

$$\sup_{P \in \mathcal{P}_j} P^n (1 - \phi_j) \leq \sqrt{\frac{\alpha}{\beta}} e^{-nj^2\varepsilon^2/9}.$$

We define $\phi$ as the supremum of all test $\phi_j$ for $j \in \mathbb{N}$. The size of this test is bounded by $\sqrt{\beta/\alpha}N(\varepsilon)\sum_{j\in\mathbb{N}}\exp[-nj^2\varepsilon^2/9]$. The power is bigger than the power of any of the tests $\phi_j$. ∎

**Lemma 6.3** (Lemma 8.1 in Ghosal et al. [2000]). *For every $\varepsilon > 0$ and probability measure $\Pi$ we have, for any $C > 0$ and $B(\varepsilon) = \{p \in \mathcal{P}: P_0 \log p_0/p \leq \varepsilon^2, P_0\big(\log p_0/p\big)^2 \leq \varepsilon^2\}$,*

$$P_0^n\Big(\int \prod_{i=1}^n \frac{p}{p_0}(X_i)\, d\Pi(P) \leq \Pi\big(B(\varepsilon)\big)e^{-(1+C)n\varepsilon^2}\Big) \leq \frac{1}{C^2 n\varepsilon^2}.$$

*Proof of Theorem 2.1.* Abbreviate $J_{n,\alpha} = n\varepsilon_{n,\alpha}^2$, so that the constant $E$ defined in the theorem is given by $E = \sup_{\alpha \gtrsim \beta_n} E_\alpha J_{n,\alpha}/J_{n,\beta_n}$.

For $\alpha \gtrsim \beta_n$ we have $B\varepsilon_{n,\beta_n} \geq B/\sqrt{H}\varepsilon_{n,\alpha} \geq \varepsilon_{n,\alpha}$. Therefore, in view of the entropy bound (2.1) and Lemma 6.2 with $\varepsilon = B\varepsilon_{n,\beta_n}$ and $\log N(\varepsilon) = E_\alpha J_{n,\alpha}$ (constant in $\varepsilon$), there exists for every $\alpha \gtrsim \beta_n$ a test $\phi_{n,\alpha}$ with,

$$P_0^n\phi_{n,\alpha} \leq \sqrt{\mu_{n,\alpha}}\frac{e^{(E_\alpha J_{n,\alpha} - KB^2 J_{n,\beta_n})}}{1 - e^{-KB^2 J_{n,\beta_n}}} \lesssim \sqrt{\mu_{n,\alpha}}e^{(E - KB^2)J_{n,\beta_n}}, \qquad (6.1)$$

$$\sup_{p\in\mathcal{P}_{n,\alpha}: d(p,p_0)\geq iB\varepsilon_{n,\beta_n}} P^n(1 - \phi_{n,\alpha}) \leq \frac{1}{\sqrt{\mu_{n,\alpha}}}e^{-KB^2 i^2 J_{n,\beta_n}}. \qquad (6.2)$$

For $\alpha < \beta_n$ we have that $\varepsilon_{n,\alpha} > \sqrt{H}\varepsilon_{n,\beta_n}$ and we cannot similarly test balls of radius proportional to $\varepsilon_{n,\beta_n}$ in $\mathcal{P}_{n,\alpha}$. However, Lemma 6.2 with $\varepsilon = B'\varepsilon_{n,\alpha}$ and $B' := B/\sqrt{H} > 1$, still gives tests $\phi_{n,\alpha}$ such that for every $i \in \mathbb{N}$,

$$P_0^n\phi_{n,\alpha} \leq \sqrt{\mu_{n,\alpha}}\frac{1}{1 - e^{-KB'^2 J_{n,\alpha}}}e^{(E_\alpha - KB'^2)J_{n,\alpha}} \lesssim \sqrt{\mu_{n,\alpha}}e^{(\underline{E} - KB'^2)J_{n,\alpha}}, \quad (6.3)$$

$$\sup_{p\in\mathcal{P}_{n,\alpha}: d(p,p_0)>iB'\varepsilon_{n,\alpha}} P^n(1 - \phi_{n,\alpha}) \leq \frac{1}{\sqrt{\mu_{n,\alpha}}}e^{-KB'^2 i^2 J_{n,\alpha}}. \qquad (6.4)$$

Let $\phi_n = \sup_{\alpha: \in A_n} \phi_{n,\alpha}$ be the supremum of all tests so constructed.

The test $\phi_n$ is more powerful than all the tests $\phi_{n,\alpha}$, and has error of the first kind $P_0^n\phi_n$ bounded by

$$P_0^n\sum_{\alpha\in A_n}\phi_{n,\alpha} \lesssim \sum_{\alpha\in A_n}\sqrt{\mu_{n,\alpha}}e^{-cJ_{n,\beta_n}},$$

for $c = (KB^2 - E) \wedge (KB^2 - \underline{E}H)$, which is bigger than 1 by assumption. Because $J_{n,\beta_n} \to \infty$ and $\sum_{\alpha\in A_n}\sqrt{\mu_{n,\alpha}} \leq \exp J_{n,\beta_n}$, this tends to zero. Consequently, for any $IB$,

$$P_0^n\Pi_n\Big(p: d(p,p_0) > IB\varepsilon_{n,\beta_n}|X_1,\ldots,X_n\Big)\phi_n \leq P_0^n\phi_n \to 0. \qquad (6.5)$$

We shall complement this with an analysis of the posterior multiplied by $1 - \phi_n$.

By Lemma 6.3 there exist events $\mathcal{A}_n$ with probability $P_0^n(\mathcal{A}_n) \geq 1 - (n\varepsilon_{n,\beta_n}^2)^{-1} \to 1$ on which

$$\int \prod_{i=1}^n \frac{p}{p_0}(X_i) \, d\Pi_n(p) \geq \lambda_{n,\beta_n} \int_{B_{n,\beta_n}(\varepsilon_{n,\beta_n})} \prod_{i=1}^n \frac{p}{p_0}(X_i) \, d\Pi_{n,\beta_n}(p)$$
$$\geq e^{-2J_{n,\beta_n}} \lambda_{n,\beta_n} \Pi_{n,\beta_n}\big(B_{n,\beta_n}(\varepsilon_{n,\beta_n})\big). \qquad (6.6)$$

Define for $i \in \mathbb{N}$,

$$S_{n,\alpha,i} = \big\{p \in \mathcal{P}_{n,\alpha} \colon iB'\varepsilon_{n,\alpha} < d(p,p_0) \leq (i+1)B'\varepsilon_{n,\alpha}\big\}, \quad \alpha < \beta_n,$$
$$S_{n,\alpha,i} = \big\{p \in \mathcal{P}_{n,\alpha} \colon iB\varepsilon_{n,\beta_n} < d(p,p_0) \leq (i+1)B\varepsilon_{n,\beta_n}\big\}, \quad \alpha \gtrsim \beta_n.$$

Then

$$\big\{p \colon d(p,p_0) > IB\varepsilon_{n,\beta_n}\big\}$$
$$\subset \bigcup_\alpha \bigcup_{i \geq I} S_{n,\alpha,i} \bigcup \bigcup_{\alpha < \beta_n} \big\{p \in \mathcal{P}_{n,\alpha} \colon IB\varepsilon_{n,\beta_n} < d(p,p_0) \leq IB'\varepsilon_{n,\alpha}\big\}$$
$$\subset \bigcup_\alpha \bigcup_{i \geq I} S_{n,\alpha,i} \bigcup \bigcup_{\alpha < \beta_n} C_{n,\alpha}\big(IB'\varepsilon_{n,\alpha}\big).$$

By Fubini's theorem and the inequality $P_0(p/p_0) \leq 1$, we have for every set $C$,

$$P_0^n \int_C \prod_{i=1}^n \frac{p}{p_0}(X_i) \, d\Pi_{n,\alpha}(p) \leq \Pi_{n,\alpha}(C),$$
$$P_0^n \int_C \prod_{i=1}^n \frac{p}{p_0}(X_i)(1-\phi_n) \, d\Pi_{n,\alpha}(p) \leq \sup_{P \in C} P^n(1-\phi_n) \, \Pi_{n,\alpha}(C).$$

Combining these inequalities with (6.6), (6.2) and (6.4) we see that,

$$P_0^n \Pi_n\big(d(p,p_0) > IB\varepsilon_{n,\beta_n} | X_1, \ldots, X_n\big)(1-\phi_n)\mathbf{1}_{\mathcal{A}_n}$$
$$\leq \sum_{\alpha \in A_n \colon \alpha \gtrsim \beta_n} \sum_{i \geq I} \frac{\lambda_{n,\alpha}}{\lambda_{n,\beta_n}} \frac{e^{-KB^2i^2 J_{n,\beta_n}} \Pi_{n,\alpha}(S_{n,\alpha,i})}{e^{-2J_{n,\beta_n}} \Pi_{n,\beta_n}\big(B_{n,\beta_n}(\varepsilon_{n,\beta_n})\big)} \frac{1}{\sqrt{\mu_{n,\alpha}}}$$
$$+ \sum_{\alpha \in A_n \colon \alpha < \beta_n} \sum_{i \geq I} \frac{\lambda_{n,\alpha}}{\lambda_{n,\beta_n}} \frac{e^{-KB'^2i^2 J_{n,\alpha}} \Pi_{n,\alpha}(S_{n,\alpha,i})}{e^{-2J_{n,\beta_n}} \Pi_{n,\beta_n}\big(B_{n,\beta_n}(\varepsilon_{n,\beta_n})\big)} \frac{1}{\sqrt{\mu_{n,\alpha}}} \qquad (6.7)$$
$$+ \sum_{\alpha \in A_n \colon \alpha < \beta_n} \frac{\lambda_{n,\alpha}}{\lambda_{n,\beta_n}} \frac{\Pi_{n,\alpha}\big(C_{n,\alpha}(IB'\varepsilon_{n,\alpha})\big)}{e^{-2J_{n,\beta}} \Pi_{n,\beta_n}\big(B_{n,\beta_n}(\varepsilon_{n,\beta_n})\big)}.$$

The third term on the right tends to zero by assumption (2.4), since $B' = B/\sqrt{H} \leq B$. We shall show that the first two terms on the right also tend to zero.

Because for $\alpha \gtrsim \beta_n$ and for $i \geq I \geq 3$ we have $S_{n,\alpha,i} \subset C_{n,\alpha}(\sqrt{2}iB\varepsilon_{n,\beta_n})$, the assumptions (2.3) shows that the first term is bounded by

$$\sum_{\alpha \in A_n: \alpha \gtrsim \beta_n} \sum_{i \geq I} \frac{\lambda_{n,\alpha}}{\lambda_{n,\beta_n}} \frac{e^{-KB^2i^2 J_{n,\beta_n}} \Pi_{n,\alpha}\big(C_{n,\alpha}(\sqrt{2}iB\varepsilon_{n,\beta_n})\big)}{e^{-2J_{n,\beta_n}}\Pi_{n,\beta_n}\big(B_{n,\beta_n}(\varepsilon_{n,\beta_n})\big)} \frac{1}{\sqrt{\mu_{n,\alpha}}}$$

$$\leq \sum_{\alpha \in A_n: \alpha \gtrsim \beta_n} \sqrt{\mu_{n,\alpha}} e^{2J_{n,\beta_n}} \sum_{i \geq I} e^{(2L-K)B^2 J_{n,\beta_n} i^2}$$

$$\leq \sum_{\alpha \in A_n: \alpha \gtrsim \beta_n} \sqrt{\mu_{n,\alpha}} e^{2J_{n,\beta_n}} \frac{e^{(2L-K)B^2 J_{n,\beta_n} I^2}}{1 - e^{(2L-K)B^2 J_{n,\beta_n}}}.$$

Because $\sum_{\alpha \in A_n} \sqrt{\mu_{n,\alpha}} \leq \exp J_{n,\beta_n}$ by assumption, this tends to zero if $(K - 2L)I^2 B^2 > 3$, which is assumed.

Similarly, for $\alpha < \beta_n$ the second term is bounded by, in view of (2.2),

$$\sum_{\alpha \in A_n: \alpha < \beta_n} \sum_{i \geq I} \frac{\lambda_{n,\alpha}}{\lambda_{n,\beta_n}} \frac{e^{-KB'^2 i^2 J_{n,\alpha}} \Pi_{n,\alpha}\big(C_{n,\alpha}(\sqrt{2}iB'\varepsilon_{n,\alpha})\big)}{e^{-2J_{n,\beta_n}}\Pi_{n,\beta_n}\big(B_{n,\beta_n}(\varepsilon_{n,\beta_n})\big)} \frac{1}{\sqrt{\mu_{n,\alpha}}}$$

$$\leq \sum_{\alpha \in A_n: \alpha < \beta_n} \sqrt{\mu_{n,\alpha}} e^{2J_{n,\beta_n}} \sum_{i \geq I} e^{(2L-K)B'^2 i^2 J_{n,\alpha}}$$

$$\leq \sum_{\alpha \in A_n: \alpha < \beta_n} \sqrt{\mu_{n,\alpha}} e^{2J_{n,\beta_n}} \frac{e^{(2L-K)B'^2 J_{n,\alpha} I^2}}{1 - e^{(2L-K)B'^2 J_{n,\alpha}}}.$$

Here $J_{n,\alpha} > HJ_{n,\beta_n}$ for every $\alpha < \beta_n$, and hence this tends to zero, because again $(K - 2L)B^2 I^2 > 3$. ∎

*Proof of Theorem 2.2.* We follow the line of argument of the proof of Theorem 2.1, the main difference being that presently $J_{n,\beta_n} = 1$ and hence does not tend to infinity. To make sure that $P_0^n \phi_n$ is small we choose the constant $B$ sufficiently large, and to make $P_0^n(\mathcal{A}_n)$ sufficiently large we apply Lemma 6.3 with $C$ a large constant instead of $C = 1$. This gives a factor $e^{-(1+C)J_{n,\beta_n}}$ instead of $e^{-2J_{n,\beta_n}}$ in the denominators of (6.7), but this is fixed for fixed $C$. The arguments then show that for an event $\mathcal{A}_n$ with probability arbitrarily close to 1 the expectation $P_0^n \Pi_n\big(d(p,p_0) > IB\varepsilon_{n,\beta_n}|X_1,\ldots,X_n\big) 1_{\mathcal{A}_n}$ can be made arbitrarily small by choosing sufficiently large $I$ and $B$. This proves the theorem. ∎

*Proof of Theorem 3.1.* The second assertion of the theorem is an immediate consequence of Theorems 2.1 and 2.2. These theorems show that the posterior concentrates all its mass on balls of radius $BI\varepsilon_{n,\beta_n}$ or $I_n\varepsilon_{n,\beta_n}$ around $p_0$, respectively. Hence the posterior cannot charge any model that do not intersect these balls.

The first assertion can be proved using exactly the proof of Theorems 2.1 and 2.2, except that the references to $\alpha \gtrsim \beta_n$ can be omitted. In the notation of the proof of Theorem 2.1 we have that

$$\bigcup_{\alpha < \beta_n} \mathcal{P}_{n,\alpha} \subset \bigcup_{\alpha < \beta_n} \bigcup_{i \geq I} S_{n,\alpha,i} \bigcup \bigcup_{\alpha < \beta_n} C_{n,\alpha}\big(IB'\varepsilon_{n,\alpha}\big).$$

It follows that $P_0\Pi_n\big(A_{n,<\beta_n}|X_1,\ldots,X_n\big)(1-\phi_n)1_{\mathcal{A}_n}$ can be bounded by the sum of the second and third terms on the right side of (6.7), which tend to zero under the conditions of Theorem 2.1, and can be made arbitrarily small under the conditions of Theorem 2.2 by choosing $B$ and/or $I$ suffcently large. ∎

## 7. Technical proofs and complements

In this section we list technical lemmas on approximation by spline spaces. We let $\|f\|_2$ and $\|f\|_\infty$ be the $L_2[0,1]$ and the supremum norm of a function $f\colon[0,1]\to\mathbb{R}$, and similarly write $\|\theta\|_2$ and $\|\theta\|_\infty$ for the Euclidean and maximum norm of $\theta\in\mathbb{R}^J$. Let $\|f\|_\alpha$ be a norm for $C^\alpha[0,1]$, for instance,

$$\|f\|_\alpha = \|f\|_\infty + \sup_{x,y\in[0,1],x\neq y}\frac{|f^{(\underline{\alpha})}(x)-f^{(\underline{\alpha})}(y)|}{|x-y|^{\alpha-\underline{\alpha}}}.$$

**Lemma 7.1.** *Let $q\geq\alpha>0$. There exists a constant $C_{q,\alpha}$ depending only on $q$ and $\alpha$ such that, for every $f$ in $C^\alpha[0,1]$,*

$$\inf_{\theta\in\mathbb{R}^J}\big\|\theta^T B_J - f\big\|_\infty \leq C_{q,\alpha}\Big(\frac{1}{J}\Big)^\alpha\|f\|_\alpha.$$

**Lemma 7.2.** *For any $\theta\in\mathbb{R}^J$,*

$$\|\theta\|_\infty \lesssim \|\theta^T B_J\|_\infty \leq \|\theta\|_\infty,$$
$$\|\theta\|_2 \lesssim \sqrt{J}\,\|\theta^T B_J\|_2 \lesssim \|\theta\|_2.$$

**Lemma 7.3.** *For any $\theta\in\mathbb{R}^J$ such that $\theta^T 1 = 0$,*

$$\underline{C_4}\|\theta\|_\infty \leq \|\log p_{J,\theta}\|_\infty \leq \overline{C}_4\|\theta\|_\infty.$$

**Lemma 7.4.** *For every $\theta_1,\theta_2\in\mathbb{R}^J$ such that $1^T(\theta_1-\theta_2)=0$,*

$$\inf_{x,\theta,J}p_{J,\theta}(x)\Big(\frac{\|\theta_1-\theta_2\|^2}{J}\wedge 1\Big) \lesssim h^2(p_{J,\theta_1},p_{J,\theta_2}) \lesssim \sup_{x,\theta,J}p_{J,\theta}(x)\Big(\frac{\|\theta_1-\theta_2\|^2}{J}\Big),$$

*where the infimum and supremum are taken over all $\theta$ on the line segment between $\theta_1$ and $\theta_2$ and all $x\in[0,1]$.*

**Lemma 7.5.** *Let $q\geq\beta$. If $\log p_0\in C^\beta[0,1]$, then the minimizer $\bar\theta_J$ of $\theta\mapsto\|\log p_{J,\theta}-\log p_0\|_\infty$ over $\theta\in\mathbb{R}^J$ with $\theta^T 1=0$ satisfies*

$$h(p_{J,\bar\theta_J},p_0) \lesssim \|\log p_{J,\bar\theta_J}-\log p_0\|_\infty \lesssim J^{-\beta}.$$

The first lemma in this list is the basic approximation lemma for splines and shows that splines of sufficient dimension are well suited to approximating smooth functions. Its proof can be found in de Boor [2001], p170. Lemmas 7.2-7.4 are (partly) implicit in Stone [1986, 1990] and can be explicitly found in Ghosal et al. [2000]. The equivalence of the $L_2$-norm or infinity-norm on the linear combinations of splines and the Euclidean or maximum norm on the coefficients (up to constants) given by Lemma 7.2 are consequences of using the B-splines, with their special properties, as a basis. Lemma 7.5 is a consequence of the other lemmas; a proof can be found in Ghosal et al. [2003].

For given $M > 0$ let $\Theta_J = \{\theta \in [-M, M]^J : \theta^T 1 = 0, \|\theta\|_\infty \leq M\}$, and write $\mathcal{P}_J$ for the set of functions $p_{J,\theta}$. By Lemma 7.3 the densities $p_{J,\theta}$ with $\theta \in \Theta_J$ take their values in the interval $[e^{-\overline{C}_4 M}, e^{\overline{C}_4 M}]$. In particular, they are uniformly bounded away from zero and infinity. Assume that the true density $p_0$ is also bounded away from 0 and infinity.

In the present case the neighbourhoods $B_{n,\alpha}(\varepsilon)$ and $C_{n,\alpha}(\varepsilon)$ defined in (1.3) take the forms $B_J(\varepsilon)$ and $C_J(\varepsilon)$ for $J = J_{n,\alpha}$ and

$$B_J(\varepsilon) = \left\{ p_{J,\theta} : P_0 \log \frac{p_0}{p_{J,\theta}} \leq \varepsilon^2, P_0 \left( \log \frac{p_0}{p_{J,\theta}} \right) \leq \varepsilon^2, \theta \in \Theta_J \right\},$$

$$C_J(\varepsilon) = \left\{ p_{J,\theta} : h(p_{J,\theta}, p_0) \leq \varepsilon, \theta \in \Theta_J \right\}.$$

Because the quotients $p_0/p_{J,\theta}$ are uniformly bounded above by $\exp(\overline{C}_4 M)$, there exists a constant $1 \leq B$, depending on $M$ only, such that

$$B_J(\varepsilon) \subset C_J(\varepsilon) \subset B_J(B\varepsilon). \tag{7.1}$$

(In fact $B$ a multiple of $M$ does; see e.g. Lemma 8 in Ghosal and van der Vaart [2007b].) In order to verify the conditions of the main theorems involving the sets $B_{n,\alpha}(\varepsilon)$ or $C_{n,\alpha}(\varepsilon)$, we can therefore restrict ourselves to the Hellinger balls $C_{n,\alpha}(\varepsilon)$. These Hellinger balls can themselves be related to Euclidean balls.

**Lemma 7.6.** *If $\theta_J$ minimizes the map $\theta \mapsto h(p_{J,\theta}, p_0)$ over $\Theta_J$ and $\varepsilon_J = h(p_0, p_{J,\theta_J})$, then there exist constants $\underline{F}$ and $\overline{F}$ such that*

$$C_J(\varepsilon) \subset \left\{ p_{J,\theta} : \theta \in \Theta_J, \underline{F} \|\theta - \theta_J\|_2 \leq \sqrt{J} 2\varepsilon \right\}, \qquad 2\varepsilon < \underline{F}, \tag{7.2}$$

$$\left\{ p_{J,\theta} : \theta \in \Theta_J, \|\theta - \theta_J\|_2 \leq \sqrt{J} \varepsilon \right\} \subset C_J(\overline{F} 2\varepsilon), \qquad \overline{F}\varepsilon \geq \varepsilon_J. \tag{7.3}$$

*Proof.* By Lemma 7.4 there exist constants $\underline{F} \leq \overline{F}$, only depending on $M$, such that, for every $\theta \in \Theta_J$,

$$\underline{F} \left( \|\theta - \theta_J\|_2 \wedge \sqrt{J} \right) \leq \sqrt{J} \, h(p_{J,\theta}, p_{J,\theta_J}) \leq \overline{F} \|\theta - \theta_J\|_2.$$

(In fact multiples of $\underline{F} = e^{-\overline{C}_4 M}$ and $\overline{F} = e^{\overline{C}_4 M}$ will do.) The set $C_J(\varepsilon)$ is empty for $\varepsilon < \varepsilon_J$. Therefore, if $p_{J,\theta} \in C_J(\varepsilon)$, then $\varepsilon \geq \varepsilon_J$ and by the triangle inequality $h(p_{J,\theta}, p_{J,\theta_J}) \leq 2\varepsilon$. If also $2\varepsilon < \underline{F}$, then the preceding display shows that $\underline{F} \|\theta - \theta_J\| \leq \sqrt{J} 2\varepsilon$. This and a similar argument for an inclusion in the other direction, yields the lemma. ∎

**Lemma 7.7.** *There exists a constant $E = E_M$ such that $\log D(\varepsilon/10, C_J(\varepsilon), h) \leq EJ$, for every $\varepsilon > 0$. (In fact $E = M \exp(\overline{C}_4 M)$ does.)*

*Proof.* If $2\varepsilon < \underline{F}$, then (7.2) shows that $C_J(\varepsilon)$ is included in the set of all $p_{J,\theta}$ with $\theta \in \overline{C}_J(\varepsilon) = \{\theta \in \Theta_J : \underline{F}\|\theta - \theta_J\|_2 \le 8\sqrt{J}2\varepsilon\}$. For given $c > 0$ there exists a constant $C$ such that we can cover $\overline{C}_J(\varepsilon)$ with $C^J$ Euclidean balls of radius $c\sqrt{J}\varepsilon$. In view of the second inequality of Lemma 7.4 this yields an $\eta$-net over $C_J(\varepsilon)$ for the Hellinger distance, for $\eta$ a multiple of $\exp(\overline{C}_4 M/2)c\varepsilon$.

If $2\varepsilon \ge \underline{F}$, then we cover $[-M, M]^J$ with of the order $(2M/c)^J$ balls of radius $c\varepsilon$ for the maximum norm. These fit in equally many Euclidean balls of radius $c\sqrt{J}\varepsilon$, and yield balls of radius a multiple of $\exp(\overline{C}_4 M/2)c\varepsilon$ in the Hellinger distance that cover $C_J(\varepsilon)$. ∎

**Lemma 7.8.** *If $\theta_J$ minimizes $\theta \mapsto h(p_0, p_{J,\theta})$ over $\theta \in \Theta_J$ and $\log p_0 \in C^\beta[0,1]$ with $\|\log p_0\|_\infty < \underline{C}_4 M$ for $\underline{C}_4$ the constant in Lemma 7.3, then $h(p_{J,\theta_J}, p_0) \lesssim J^{-\beta}$.*

*Proof.* In view of Lemma 7.5 it suffices to show that $\bar{\theta}_J$ defined there satisfies $\|\bar{\theta}_J\|_\infty \le M$. By the triangle inequality, $\underline{C}_4\|\bar{\theta}_J\|_\infty \le \|\log p_{J,\bar{\theta}_J} - \log p_0\|_\infty + \|\log p_0\|_\infty$, where the first term on the right is of order $O(J^{-\beta})$ by Lemma 7.5. ∎

**Lemma 7.9.** *If $v_J$ is the volume of the $J$-dimensional unit ball, then $J \mapsto \sqrt{J}^J v_J$ is increasing, and, as $J \to \infty$,*

$$\sqrt{J}^J v_J = \frac{\sqrt{J}^J \sqrt{\pi}^J}{\Gamma(J/2 + 1)} = \frac{\sqrt{2\pi e}^J}{\sqrt{\pi J}}(1 + o(1)).$$